\documentclass[12pt]{amsart}

\usepackage[utf8]{inputenc}
\usepackage[T1]{fontenc}

\usepackage[a4paper,margin=2.5cm]{geometry}
\usepackage{amssymb,amsfonts,amstext}
\usepackage{thmtools}
\usepackage{mathtools}
\usepackage{mathrsfs}
\usepackage{mathabx}
\usepackage{scalerel}
\usepackage{relsize}
\usepackage{faktor}
\usepackage{xfrac}
\usepackage{blkarray}

\usepackage{graphicx}
\usepackage{float}
\usepackage[tableposition=top]{caption}
\usepackage{subfig}

\usepackage{array}
\usepackage{multirow}
\usepackage{enumerate}
\usepackage{enumitem}

\newcolumntype{L}[1]{>{\raggedright\let\newline\\\arraybackslash\hspace{0pt}}m{#1}}
\newcolumntype{C}[1]{>{\centering\let\newline\\\arraybackslash\hspace{0pt}}m{#1}}
\newcolumntype{R}[1]{>{\raggedleft\let\newline\\\arraybackslash\hspace{0pt}}m{#1}}

\makeatletter
\renewcommand\section{%
  \@startsection{section}{1}\z@
  {.7\linespacing \@plus \linespacing}
  {.5\linespacing}{\normalfont\bfseries\scshape\centering}%
}

\renewcommand\subsubsection{%
\@startsection{subsubsection}{3}\z@ {.5\linespacing \@plus .7\linespacing }{-.5em}{\normalfont\bfseries}
}

\usepackage[most]{tcolorbox}

\newtcbox{\catboxI}{
  on line,
  colback=white,
  colframe=black,
  boxrule=0.5pt,
  arc=4pt,
  left=4mm,
  right=4mm,
  top=2mm,
  bottom=2mm,
  width=1cm,
  boxsep=0pt,
  valign=center,
  halign=center
}
\newtcbox{\catboxII}{
  on line,
  colback=gray!5,
  colframe=black,
  boxrule=0.5pt,
  arc=2pt,
  left=6mm,
  right=6mm,
  top=2mm,
  bottom=2mm,
  width=1cm,
  boxsep=0pt,
  valign=center,
  halign=center
}
\newtcbox{\catboxIII}{
  on line,
  colback=gray!5,
  colframe=black,
  boxrule=0.6pt,
  arc=3pt,
  left=2mm,
  right=2mm,
  top=2mm,
  bottom=2mm,
  width=1cm,
  boxsep=0pt,
  valign=center
}

\usepackage{tikz}
\usepackage{tikz-cd}
\usetikzlibrary{arrows,arrows.meta,decorations.markings,positioning}

\usepackage[nospace,noadjust]{cite}
\usepackage[colorlinks=true,
            linkcolor=green!30!black,
            citecolor=blue,
            urlcolor=blue,
            pagebackref=true]{hyperref}
\usepackage{aliascnt}
\usepackage{cleveref}

\declaretheoremstyle[
  headfont=\bfseries,
  bodyfont=\itshape,
]{plainstyle}

\declaretheoremstyle[
  headfont=\bfseries,
  bodyfont=\normalfont,
]{defstyle}

\declaretheorem[
  style=plainstyle,
  numberwithin=section,
  name=Theorem
]{thm}

\declaretheorem[sibling=thm, style=defstyle, name=Remark]{remark}

\declaretheorem[
numbered=no,style=defstyle, name=Question]{question}

\declaretheorem[name=Claim,style=plainstyle]{claim}

\declaretheorem[numbered=no, name=Claim]{claim*}
\declaretheorem[numbered=no, name=Problem]{problem*}

\declaretheorem[style=plainstyle, name=Conjecture]{conjecture}

\crefname{thm}{theorem}{theorems}
\Crefname{thm}{Theorem}{Theorems}

\crefname{prop}{proposition}{propositions}
\Crefname{prop}{Proposition}{Propositions}

\crefname{lemma}{lemma}{lemmas}
\Crefname{lemma}{Lemma}{Lemmas}

\crefname{fact}{fact}{facts}
\Crefname{fact}{Fact}{Facts}

\crefname{corollary}{corollary}{corollaries}
\Crefname{corollary}{Corollary}{Corollaries}

\crefname{conjecture}{conjecture}{conjectures}
\Crefname{conjecture}{Conjecture}{Conjectures}

\crefname{definition}{definition}{definitions}
\Crefname{definition}{Definition}{Definitions}

\crefname{remark}{remark}{remarks}
\Crefname{remark}{Remark}{Remarks}

\crefname{example}{example}{examples}
\Crefname{example}{Example}{Examples}

\crefname{exercise}{exercise}{exercises}
\Crefname{exercise}{Exercise}{Exercises}

\crefname{question}{question}{questions}
\Crefname{question}{Question}{Questions}

\crefname{claim}{claim}{claims}
\Crefname{claim}{Claim}{Claims}

\crefname{case}{case}{cases}
\Crefname{case}{Case}{Cases}

\crefname{problem}{problem}{problems}
\Crefname{problem}{Problem}{Problems}

\usepackage{hebcal}

\newcommand{\interior}[1]{{\kern0pt#1}^{\mathrm{o}}}
\newcommand*{\medcap}{\mathbin{\scalebox{1.5}{\ensuremath{\cap}}}}

\newcommand\preceqdot{\mathrel{\mathcal{O}align{$\preceq$\cr
  \hidewidth\raise0.125ex\hbox{$\cdot\mkern0.5mu$}\cr}}}
\newcommand\precdott{\mathrel{\mathcal{O}align{$\prec$\cr
  \hidewidth\raise0.015ex\hbox{$\cdot\mkern0.5mu$}\cr}}}

\newcommand\Item[1][]{%
  \ifx\relax#1\relax \item \else \item[#1] \fi
  \abovedisplayskip=0pt\abovedisplayshortskip=0pt~\vspace*{-\baselineskip}
}

\def\zerovec{{\bf 0}}

\def\b{{\bf b}}
\def\x{{\bf x}}
\def\y{{\bf y}}

\def\p{\boldsymbol{p}}
\def\m{{\bf m}}
\def\n{{\bf n}}
\def\u{{\bf u}}
\def\v{{\bf v}}
\def\N{{\bf N}}

\def\R{{\mathbb{R}}}
\def\Z{{\mathbb{Z}}}
\def\bbQ{\mathbb{Q}}

\def\cA{\mathcal{A}}
\def\cL{\mathcal{L}}
\def\cW{\mathcal{W}}
\def\cP{\mathcal{P}}

\def\cO{{\mathcal{O}}}

\def\vol{{\rm vol}}

\def\vis{{\rm vis}}

\def\SL{{\rm SL}}
\def\ASL{{\rm ASL}}
\def\Cl{\mathscr{C}}
\def\scrX{\mathscr{L}}
\def\scrY{\mathscr{A}}

\newcommand{\rmd}{\mathrm{d}}

\newcommand{\df}{{\, \stackrel{\mathrm{def}}{=}\, }}
\def\cW{\mathcal{W}}
\def\cL{\mathcal{L}}

\newcommand{\card}[1]{\#\left(#1\right)}
\renewcommand{\exp}[1]{\mathrm{exp}\left(#1\right)}

\DeclareMathAlphabet{\mathpzc}{OT1}{pzc}{m}{it}
\newcommand{\mink}{\mathpzc{m}}

\usepackage{xcolor}
\usepackage[textsize=tiny]{todonotes}
\makeatletter
\def\keywords{\xdef\@thefnmark{}\@footnotetext}
\makeatother

\title[Diffraction spectrum of visible points]{On the diffraction spectrum of the set of visible points in lattices and certain cut-and-project sets}

\author[R. Kumar]{Rishi Kumar}
\address{School of Mathematical Sciences\\
Tel Aviv University\\
Tel Aviv\\
69978\\
Israel}
\email{rkumar@tauex.tau.ac.il}

\author[C. Ospina]{Carlos Ospina}
\address{School of Mathematical Sciences\\
Tel Aviv University\\
Tel Aviv\\
69978\\
Israel}
\email{ospina.math@icloud.com}

\date{}

\begin{document}

\begin{abstract}
Let $k\geq 2$ be a positive integer. It is known that the set of visible lattice points from the origin in $\mathbb{Z}^k$ has a translation bounded pure point diffraction spectrum. We investigate these properties for  sets of points simultaneously visible from a finite set of lattice points $ \{\x_1,\dots,\x_n\} \subseteq \mathbb{Z}^k$. We provide explicit formulas for the coefficients of the diffraction spectrum. Additionally, we generalize our procedure to show that the set of visible points from the origin in certain classes of cut-and-project sets has a  translation bounded pure point diffraction spectrum.
\end{abstract}

\maketitle

\section{Introduction}
Let $\Lambda \subset G$ be a uniformly discrete subset of a locally compact $\sigma$-compact Abelian group. One approach to studying the structure of $\Lambda$ is through the analysis of the \emph{autocorrelation} and  \emph{diffraction spectrum} associated with the \emph{Dirac comb} measure
\[
\omega_{\Lambda} \df \sum_{\x \in \Lambda} \delta_\x,
\]
which is a translation bounded measure (see \cite[Example~8.6]{Baake1} in the Euclidean case and \cite[Lemma 2]{BAAKE_LENZ_2004} for a general locally compact $\sigma$-compact Abelian group $G$). We recall these definitions in \Cref{sec:diff_spectrum}.

To define the autocorrelation and the diffraction spectrum, one must make several choices, including an increasing sequence $\mathcal{A}= \{A_n\}_{n\in \mathbb{N}}$ of measurable subsets of $G$ and a normalized Haar measure $\theta_G$. 
We denote by $\omega_\Lambda|_{A_n}$ the restriction of $\omega_\Lambda$ to the set $A_n$ and denote by $\widetilde{\omega_\Lambda}$ the twisted version of $\omega_\Lambda$; see \cref{eq:twist_measure} for the definition.
An \emph{autocorrelation} is a weak-$*$ accumulation point of the family of convolutions
\[
\left\{
\frac{
\omega_\Lambda|_{A_n}
\bigast
\widetilde{\omega_\Lambda|_{A_n}}
}{
\theta_G(A_n)
}
:
n\in \mathbb{N}
\right\}
\]
as $n \to \infty$. When this family has a unique accumulation point, that measure is called the \emph{autocorrelation} of $\omega_\Lambda$ and is denoted by $\gamma_{\omega_\Lambda}$. The precise choices of these objects will be specified later in the context of cut-and-project sets.
The \emph{diffraction spectrum} of $\Lambda$ is the Fourier transform of the autocorrelation,
\[
\widehat{\gamma_{\omega_\Lambda}}.
\]
If the support of the diffraction spectrum is discrete, we say that the diffraction spectrum is \emph{pure point}.

For more details on diffraction spectra, see \cite{Huck1,Martin,Moody,Hof2} and the references therein.

\subsection{Results}
Let $\Lambda \subset \mathbb{R}^d$ be uniformly discrete point set, and let $D\subset \mathbb{R}^d$ be a bounded measurable set with $\vol(D)>0$. Denote $TD= \{t\x\, :\, \x \in D,\, t\in [0,T],\, T>0\}$ for $T>0$ to be the dilated set. The density of $\Lambda$ with respect to $D$ is defined to be
$$\theta(\Lambda)\df \lim_{t\to \infty}\frac{\#(\Lambda\medcap TD)}{\vol(TD)},$$
provided the limit exists, where $\card{A}$ is the cardinality of a set~$A$. See \eqref{def:density_of_a_set} for a more general definition of point sets in $\sigma$-compact locally compact abelian groups.  We denote
$$\Lambda_{\vis}=\{\y \in \Lambda\setminus \{\textbf{0}\}\, :\, t\y \not \in \Lambda,\, \forall t\in (0,1)\},$$
the set of points visible from the origin. 

We now describe the main results of this paper. We begin with the classical case of visible points in the integer lattice $\mathbb{Z}^k$, which motivates our study of points that are simultaneously visible from a finite set of lattice points and, later, visible points from the origin in some class of cut-and-project sets.

\subsubsection{Visible points in the integer lattice}

Let $k\geq 2$ be a positive integer. Denote the  Dirac comb of $\mathbb{Z}^k$ by
\[
\omega_{\mathbb{Z}^k} = \sum_{\x \in \mathbb{Z}^k} \delta_{\x}.
\]
By the Poisson summation formula, the diffraction spectrum is pure point and is given by
\[
\widehat{\gamma_{\omega_{\mathbb{Z}^k}}} = \omega_{\mathbb{Z}^k}.
\]
In particular, the coefficients of the diffraction spectrum are constant and equal to $1$. Physically, if X-rays (photons of a fixed wavelength) are transmitted through $\mathbb{Z}^k$, one observes a periodic array of bright spots known as \emph{Bragg peaks}. These Bragg peaks are supported on $\mathbb{Z}^k$, and each peak has unit intensity.

Given two distinct points $\x = (x_1,\ldots,x_k)$ and $\y = (y_1,\ldots,y_k)$ in $\mathbb{Z}^k$, we say that $\y$ is \emph{mutually visible} from $\x$ if no other point of $\mathbb{Z}^k$ lies on the line segment joining them. 
This occurs if and only if $\gcd(x_1- y_1, \ldots, x_k -y_k) =1$ where $\gcd$ stands for greatest common divisor, see \cite{Rea}.
Let $V_k= \mathbb{Z}^k_{\vis}$ denote the set of points in $\mathbb{Z}^k$ that are visible from the origin, namely,
\[
V_k = \{(x_1,\ldots,x_k)\in \mathbb{Z}^k\,:\, \gcd(x_1,\ldots,x_k)=1\}.
\]

By the Chinese Remainder Theorem, the set $\Z^k$ contains arbitrarily large cubes containing no visible points (cf. \cite[Theorem~5.29]{Apostol} for the integer lattice, and \cite[Proposition~4]{Moody} for arbitrary full-rank lattices). Thus, despite being uniformly discrete, $V_k$ fails to be relatively dense and hence is not a Delone set (see \Cref{subsec:Generalized_cut_and_project_sets}). In fact, since these holes repeat periodically throughout the lattice, adjoining any set of zero density to $V_k$ still does not produce a Delone set.

\begin{figure}[htbp]
\centering
\includegraphics[width=0.50\linewidth]{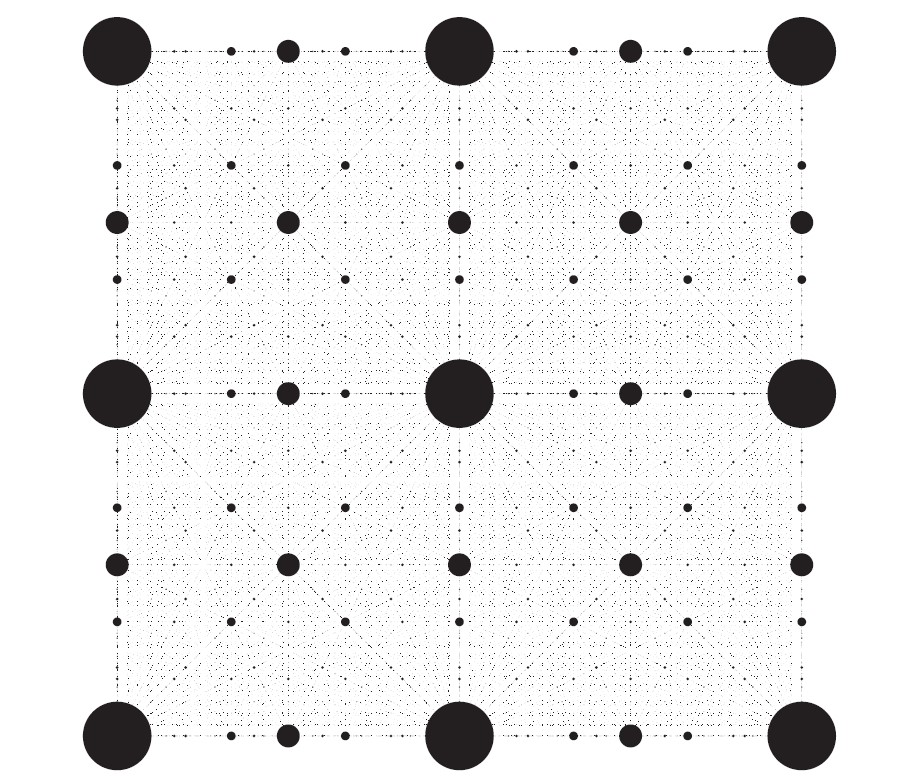}
\caption{Diffraction of the visible points of $\mathbb{Z}^2$ \cite[p.~424]{Baake1}.}
\label{fig:visible_diffraction}
\end{figure}

For $k \geq 2$,  Baake--Moody--Pleasants \cite[Theorem~3]{Moody} studied the diffraction spectrum $\widehat{\gamma_{\omega_{V_k}}}$ of the set $V_k \subset \mathbb{Z}^k$. They showed that $\widehat{\gamma_{\omega_{V_k}}}$ is again a translation bounded pure point measure. Moreover, the diffraction spectrum can be computed explicitly and is supported on the subset of $\mathbb{Q}^k$ consisting of vectors with square-free denominators and is given by
\[
\widehat{\gamma_{\omega_{V_k}}}= \sum_{\substack{\u=\frac{\m}{N}\in \mathbb{Q}^k\\ \m \in \Z^k \\N\, \mathrm{square\text{-}free}
\\\mathrm{gcd}(N,\m)=1
}}|a(\textbf{u})|^2\delta_{\textbf{u}}, \qquad a(\textbf{u})= \frac{1}{\zeta(k)}\prod_{p \mid N}\frac{1}{1-p^k}.
\]
In fact, they established this result in the more general setting of arbitrary full-rank lattices in dimension $k\geq 2$, as well as for the one-dimensional set of $k$th-power-free integers.
See \Cref{fig:visible_diffraction} for an illustration of the diffraction spectrum of visible lattice points in $\mathbb{Z}^2$.
The study of diffraction spectra for visible lattice points and $k$-free lattice points has been further extended to $\mathcal{B}$-free lattice systems (see, for example \cite[\S~4]{Luz}, \cite[Corollary 21]{Strungaru1} and references therein).
This motivates the following question. Let $S \subset \mathbb{Z}^k$ be a finite set.
\begin{question}
What can be said about the diffraction spectrum of the points in $\mathbb{Z}^k$ that are simultaneously visible from every element of $S$?
\end{question}
Denote by $\mathcal{P}$ the set of prime numbers in $\mathbb{N}$.
For each $p \in \mathcal{P}$, the map $\pi_p$ denotes the natural projection
\begin{equation}\label{eq:projection_mod_p}
\pi_p : \Z^k \longrightarrow (\Z/p\Z)^k, \qquad
(x_1, \dots, x_k) \longmapsto \bigl(x_1 \bmod{p}, \dots, x_k \bmod{p}\bigr),
\end{equation}
and
\[
s(p) \df \card{\pi_p(S)}
\]
is the cardinality of the image of $S$ under $\pi_p$. 
Define $V(S) \subseteq \mathbb{Z}^k$ to be the set of integer lattice points that are simultaneously visible from all points of $S$.  Rumsey~\cite[Theorem~7]{Rum} proved that the density of the points simultaneously visible from $S$, averaged over the boxes $[-n,n]^k$, is given by
\begin{equation}\label{eq:sp}
\theta(V(S)) = \prod_{p \in \mathcal{P}} \left(1 - \frac{s(p)}{p^k}\right).
\end{equation}
 In the special case when points of $S$ are visible to each other, the density result \eqref{eq:sp} was proved earlier by Rearick \cite{rearick1966mutually} with an error term.
Rumsey extended this formula to certain infinite subsets $S$ satisfying appropriate conditions and observed that the same density formula holds when averaging over rectangular boxes
\[
[-n_1, n_1] \times \cdots \times [-n_k, n_k],
\]
as $\min\{n_1, \ldots, n_k\} \to \infty$. The density results when averaging over Euclidean balls follow from \cite[Theorem 1.3]{Takeda} and  \cite[Lemma 1]{pleasants2013entropy}. We refer the reader to \cite{Berend2} and references therein for more details on simultaneous visibility.
Visible lattice points have appeared in several different contexts. Besides diffraction and dynamical questions, geometric and probabilistic properties of visible points in high-dimensional cubes were investigated recently by Athreya, Cobeli, and Zaharescu in \cite{AthreyaCobeliZaharescu2023}.
They studied visibility phenomena in high-dimensional cubes and showed that almost all self-visible triangles with vertices in $[0,N]^d$ are nearly equilateral, with side lengths close to $\frac{N\sqrt{d}}{\sqrt{6}}$.

In the present work, we consider the Dirac comb measure
\[
\omega_{V(S)} \df \sum_{\x \in V(S)}\delta_\x
\]
and prove the following.

\begin{thm}\label{thm:DS_VS}
Let $S \subset \mathbb{Z}^k$ be a finite set with $k \geq 2$ and assume that the density of $V(S)$, computed along the averaging sets $[-L,L]^k$, is strictly positive. Then
\begin{enumerate}
    \item\label{part:1:thm:DS_VS}
    The diffraction spectrum $\widehat{\gamma_{\omega_{V(S)}}}$ of the set of visible points $V(S)$ is a translation bounded pure point measure;

    \item\label{part:2:thm:DS_VS}
    More precisely,
    \[
    \widehat{\gamma_{\omega_{V(S)}}}
    =
    \sum_{\substack{
    \u=\frac{\m}{N}\in \mathbb{Q}^k \\
    \m \in \mathbb{Z}^k \\
    N \in \mathbb{N}\ \mathrm{square\text{-}free
    }\\
    \mathrm{gcd}(N,\m)=1}}
    |a(\u)|^2\delta_{\u},
    \]
    where
    \[
    a(\u)=
    \theta(V(S))
    \prod_{\substack{p\in\mathcal{P}\\ p\mid N}}
    \frac{
    \left(
    -\frac{1}{p^k}
    \sum_{\x\in \pi_p(S)}
    \exp{- 2\pi i\frac{\langle \v_p,\x\rangle}{p}}
    \right)
    }{
    \left(
    1-\frac{s(p)}{p^k}
    \right)
    },\qquad \v_{p}= -\m \bmod{p}.
    \]
\end{enumerate}
\end{thm}

The sets $V(S)$ considered in  \Cref{thm:DS_VS} belongs to a broader class of examples arising from generalized Erd{\"o}s sieves.   These sieves were studied recently in \cite{araujo2026sarnak} and include the $\mathcal{B}$-free integers as a special case; see also \cite[\S~2]{Glundbach2024} and the references therein.

We prove \Cref{thm:DS_VS} in \Cref{proof:thm:DS_VS,proof:thm:DS_VS_2}.

\subsubsection{Visible points in cut-and-project sets}
Let $n = d+m$ for positive integers $n,d,m$, and let
$\pi_{\textup{phys}}$ and $\pi_{\textup{int}}$ denote the
projections onto the first $d$ coordinates and the last $m$ coordinates, respectively, in the direct sum decomposition
\begin{equation}\label{eq: direct sum decomposition1}
  \R^n= \R^d\oplus \R^m.
\end{equation}
Let $\cW\subset \mathbb{R}^m$ be a compact measurable subset referred to as a window. The cut-and-project set $\Lambda(\cW,\cL)$ is the projection onto $\R^d$ of the lattice points that belong to the strip $\R^d \times \cW$, that is,
\[
\Lambda(\cW,\cL) 
\df 
\{\x \in \R^d : (\x,\y) \in \cL,~\y \in \cW\}.
\]
We say that $\Lambda=\Lambda(\cW,\cL)$ is \emph{irreducible} if $\pi_{{\rm int}}(\cL)$ is dense in $\mathbb{R}^m$, $\pi_{{\rm phys}}|_{\cL}$ is injective and the boundary of the window has volume zero.
See \Cref{subsec:Generalized_cut_and_project_sets} for a more general definition of a cut-and-project set $\Lambda \subset G$ for $\sigma$-compact locally compact abelian groups.

\begin{figure}[htbp]
\centering
{\includegraphics[width=0.35\linewidth]{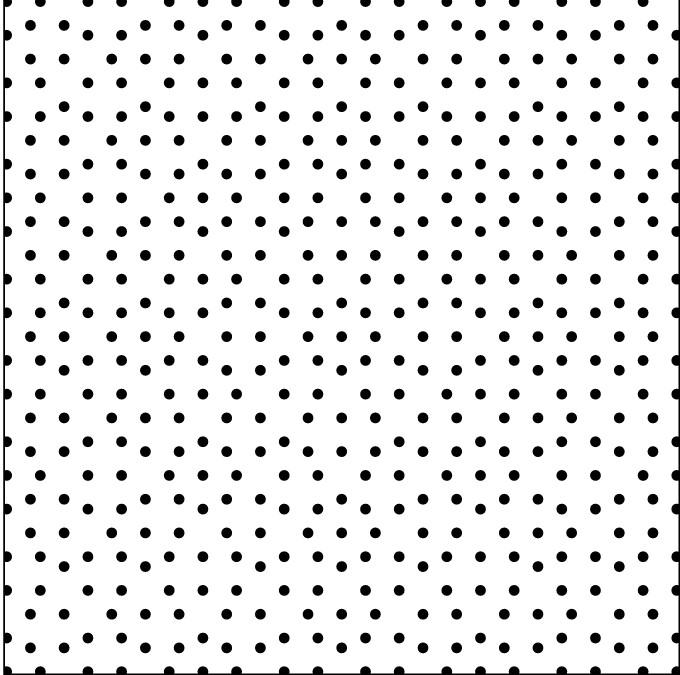}}\qquad
{\includegraphics[width=0.41\linewidth]{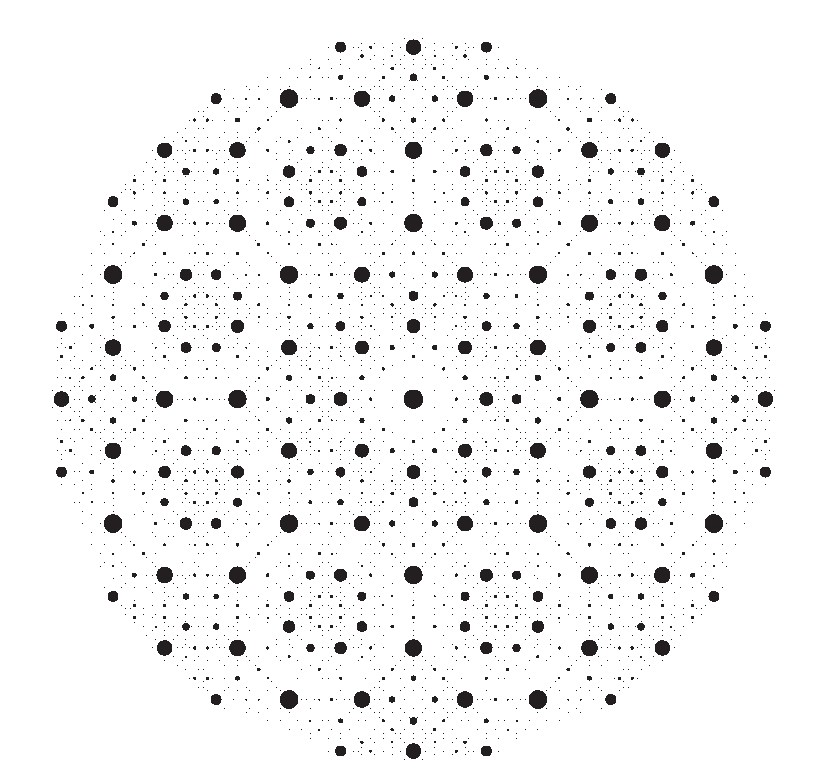}}\\
\caption{Ammann--Beenker point set (left) and its diffraction image (right), adapted from \cite{Sing_slides} and \cite[p.~372]{Baake1}, respectively.}
\label{fig15}
\end{figure}

The Dirac comb associated with the cut-and-project set $\Lambda \df \Lambda(\cW,\cL)$ is defined by
\[
\omega_{\Lambda} \df \sum_{\x \in \Lambda} \delta_{\x}.
\]
The diffraction spectrum $\widehat{\gamma_{\omega_{\Lambda}}}$ of a regular cut-and-project set is a translation bounded pure point measure, provided that $\mathcal{W}$ is compact. See \cite[Theorem~9.4]{Baake1} and \cite{Hof2} for the Euclidean case. For the general setting of locally compact Abelian groups, see \cite{Martin,BaakeMoody2}. Moreover, it admits an explicit description in terms of the Fourier transform of the window. See \Cref{fig15} for the Ammann--Beenker point set and its diffraction spectrum.

Recall that the visible points of $\Lambda(\cW,\cL)$ are precisely those points whose segment to the origin contains no other point of $\Lambda(\cW,\cL)$. Equivalently,
\[
\Lambda_{\vis}= \Lambda(\cW,\cL)_{\vis}\df \{\x\in \Lambda(\cW,\cL)\setminus \{\zerovec\}\, :\, t\x\not\in \Lambda(\cW,\cL),\, \forall t\in (0,1)\}.
\]
Marklof and Str{\"o}mbergsson \cite[Theorem~1]{Marklof3} proved that for any regular cut-and-project set, the density of $\Lambda_{\vis}$ exists and satisfies
\begin{equation}\label{eq:MS_inequalities}
0 < \theta(\Lambda_{\vis}) \leq \theta(\Lambda).
\end{equation}

On the other hand, Hammarhjelm \cite{Gustav} constructed examples of cut-and-project sets $\Lambda(\cW,\cL) \subseteq \R^2$ for which the inequality in \eqref{eq:MS_inequalities} is strict and computed $\theta(\Lambda_{\vis})$ explicitly. The examples considered by Hammarhjelm included the Amman-Beenker point set and sets associated with the Penrose tiling vertex set. Some of these examples were also considered by Sing \cite{Sing_slides}. The examples considered by Hammarhjelm and Sing arise from algebraic lattices over real quadratic fields.
To state our results, we recall some definitions from algebraic number theory.\\
Let $K = \bbQ(\sqrt{\mathfrak{d}})$ be a real quadratic field with square-free $\mathfrak{d} \geq 2$, and let $\mathcal{O}_K$ be its ring of integers. We assume throughout that $\mathcal{O}_K$ is a principal ideal domain. It is conjectured that there are infinitely many real quadratic fields $K$ for which $\mathcal{O}_K$ is a principal ideal domain (see \cite[p.~37]{Neukirch1}).\\
For an ideal $I \subseteq \cO_K$, the \emph{norm} of $I$ is defined by
$$
N(I) \df \card{\cO_K/I},
$$
which is always finite.
We denote the Dedekind zeta function by
$$
\zeta_{\cO_K}(s) \df \sum_{I\subseteq \cO_K}\frac{1}{N(I)^s},
\qquad
\Re(s)>1.
$$\\
The Minkowski embedding
$$
f: K\hookrightarrow \mathbb{R}^2, \qquad \alpha \mapsto (\alpha, \sigma(\alpha)),
$$
embeds $\cO_K$ as a lattice $f(\cO_K)$ in $\mathbb{R}^2$, where $\sigma$ is the non-trivial Galois automorphism of $K$. Let $\lambda>1$ denote the fundamental unit of $\cO_K$, equivalently, the generator of the free part of $\cO_K^\times$. We say that $K$ satisfies the \emph{Hammarhjelm condition} if
$$
f(\mathcal{O}_K)\medcap ((1,\lambda)\times [-1,1])= \varnothing.
$$\\
For $2\leq \mathfrak{d}\leq 100$, the Hammarhjelm condition is satisfied only for $\mathfrak{d} = 2, 5, 13, 29, 53$; see \cite{Barak2}. In private communication with the authors, Z. Rudnick indicated that he proved that there are only finitely many real quadratic fields $K$ for which the Hammarhjelm condition holds.

More recently, in \cite{Barak2}, the first-named author together with I. Gringlaz and B. Weiss extended Hammarhjelm’s density results to all dimensions $d \geq 2$ and established an error term.  Moreover, \cite[Theorem~1.1]{Barak2}  
implies that $\Lambda_{\vis}$ contains arbitrarily large holes. In fact, it follows from the proof of \Cref{thm:DS_of_visible_points_in_CPS} together with \cite[Proposition~2.12]{Huck1} that these holes repeat. Consequently, adjoining any set of zero density to $\Lambda_{\vis}$ still does not produce a Delone set.

\begin{thm}\emph{\cite[Theorem~5.1]{Barak2}}\label{theorem on error term}
Let $d \geq 2$, let $m=d$ and let $\mathcal{W}\subset \R^{d}$ be a convex centrally symmetric window. Let $K$ be a real quadratic field with ring of integers $\mathcal{O}_K$. Assume that $K$ satisfies the Hammarhjelm condition and
\begin{equation}\label{eq:def_lattice}
\mathcal{L} \df \{(x_1,\ldots, x_d,\sigma(x_1), \ldots, \sigma(x_d)) :
x_i \in \mathcal{O}_K\} \subset \R^{2d}.
\end{equation}

Let $\Lambda = \Lambda(\mathcal{W}, \mathcal{L})$ and let $D \subset \R^d$ be a convex averaging set. Then, as $T \to \infty$,
\begin{equation*}
\begin{split}
\frac{\card{\Lambda_{\vis} \medcap TD}}{\vol(TD)}
=
\left(1-\frac{1}{\lambda^d}\right)
\cdot
\frac{\theta(\Lambda)}{\zeta_{\mathcal{O}_K}(d)}
+
\begin{cases} 
O\left(\frac{\log T}{\sqrt{T}}\right),& \qquad d=2,\\
O\left(\frac{1}{\sqrt{T}}\right),& \qquad d\geq 3.
\end{cases}
\end{split}
\end{equation*}
\end{thm}

We prove the following.

\begin{thm}\label{thm:DS_of_visible_points_in_CPS}
Let $K$ be a real quadratic field satisfying the Hammarhjelm condition. Let $\mathcal{W} \subset \R^d$ be a convex centrally symmetric window, and let $\mathcal{L} \subset \R^{2d}$ be the lattice defined in \eqref{eq:def_lattice}. Let $\Lambda = \Lambda(\mathcal{W}, \mathcal{L})$ be the associated cut-and-project set and let $D\subset \mathbb{R}^d$ be a convex Jordan measurable averaging set. Then the diffraction spectrum $\widehat{\gamma_{\omega_{\Lambda_{\vis}}}}$ is a translation bounded pure point measure along the sequence $\{TD\}_{T >0}$.
\end{thm}

We prove the statement in \Cref{proof:thm:DS_of_visible_points_in_CPS}. The main ingredients are the density estimates from \cite{Barak2} and Baake, Huck, and Strungaru's results on the diffraction of weak model sets.

\subsubsection{Random cut-and-project sets}
Let $\SL_n(\R)$ be the group of linear transformations of $\R^n$ that are volume-preserving and orientation preserving. Let $\ASL_n(\R) = \SL_n(\R) \ltimes \R^n $ be the group of affine transformations on $\R^n$ that are volume-preserving and orientation preserving. We consider the spaces
\[
\scrX_n \df \SL_n(\R)/\SL_n(\Z) \quad \text{and} \quad \scrY_n \df \ASL_n(\R)/\ASL_n(\Z)
\]
with the quotient topology, which parameterize covolume-one lattices and grids in $\R^n$, respectively. Let $m_{\scrX_n}$ and $m_{\scrY_n}$ be the Haar-Siegel measures on $\scrX_n$ and $\scrY_n$, namely, the unique probability measures invariant under the actions of $\SL_n(\R)$ and $\ASL_n(\R)$, respectively.

Consider the collection $\Cl(\R^n)$ of closed subsets of $\R^n$. Under the \emph{ Chabauty-Fell topology}, $\Cl(\R^n)$ is compact, see \cite[\S~2.2]{Barak1} and the references therein.
The restriction of the Chabauty-Fell topology to the sets $\scrX_n$ and $\scrY_n$ coincides with the quotient topology.
This gives the natural embeddings $\scrX_n \subset \scrY_n \subset \Cl(\R^n)$.

Fixing the direct sum \eqref{eq: direct sum decomposition1}, its corresponding projections $\pi_{\mathrm{phys}}$ and $\pi_{\mathrm{int}}$, and a window $\mathcal{W} \subset \R^m$, define the map
\[
\Psi: \scrY_n \to \Cl(\R^d), \qquad \Psi(\mathcal{L}) \df \Lambda(\mathcal{W}, \mathcal{L}).
\]
Let $\bar{\mu}$ and $\mu$ denote the pushforwards of the measures $m_{\scrY_n}$ and $m_{\scrX_n}$ by the map $\Psi$, respectively. Both measures $\bar{\mu}$ and $\mu$ define random cut-and-project sets by fixing the decomposition \eqref{eq: direct sum decomposition1} and the window $\cW$, and then selecting a grid or lattice $\mathcal{L}$ according to the corresponding Haar-Siegel measure.

Since the groups $\ASL_d(\R)$ and $\SL_d(\R)$ act on $\R^d$, the action extends naturally to $\Cl(\R^d)$. A probability measure on $\Cl(\R^d)$ is a \emph{Ratner-Marklof-Str\"ombergsson measure} (RMS measure for short) if it is invariant and ergodic with respect to the $\SL_d(\R)$-action and gives full measure to the set of irreducible cut-and-project sets. Such a measure is called an \emph{RMS linear measure}. Since $\SL_d(\R)$ can be seen as a subgroup of $\ASL_d(\R)$, an RMS linear measure does not need to be $\ASL_d(\R)$-invariant. If it is $\ASL_d(\R)$-invariant, we call it an \emph{RMS affine measure}. The RMS measures are completely classified; see \cite{Marklof2, Barak1}.

In \cite[Propositions~3.1~and~3.2]{Barak2}, it is shown that for $\bar{\mu}$-almost every cut-and-project set, we have the equality
\[
\Lambda(\mathcal{W},\mathcal{L}) \setminus \{\textbf{0}\} = \Lambda(\mathcal{W},\mathcal{L})_{\vis}.
\]
If $\cL$ is an affine lattice, i.e. $\cL=\cL_0+ \x$ for some lattice $\cL_0\subset \mathbb{R}^n$ and $\x\in \mathbb{R}^n$, then
\[
\Lambda(\cW,\cL)= \Lambda(\cW-\pi_{\rm int}(\x),\cL_0)+ \pi_{\rm phys}(\x).
\]
Recall that the map $\Psi$ depends on the decomposition on $\eqref{eq: direct sum decomposition1}$ and the compact window $\cW$ with nonempty interior and boundary of Lebesgue measure zero. Hence, it follows from \cite[Theorem~9.4]{Baake1} and \cite{Hof2} that the diffraction spectrum of $\bar{\mu}$-almost every cut-and-project set is a translation bounded pure point measure along $B_R(0)$ and $[-L,L]^d$, respectively.

Let $d\geq 2$, $\mathcal{W}\subset \mathbb{R}^m$ be star-shaped with respect to the origin, and $\mathcal{A}=\{A_t\}_{t>0}$ be an unbounded ordered family in $\R^d$ (see \Cref{subsec:van_Hove_sequences} for the definition). It follows from \cite[Theorem~3.4]{Barak2} that the density along $\mathcal{A}$ of $\mu$-almost every cut-and-project set exists and is given by
\begin{equation}\label{eq. density of the random cut-and-project set}
\theta(\Lambda(\mathcal{W},\mathcal{L}))= \frac{\vol(\mathcal{W})}{\zeta(n)}.
\end{equation}
In \Cref{subsec:van_Hove_sequences}, we show that every van Hove sequence can be completed to an unbounded ordered family. Consequently, \eqref{eq. density of the random cut-and-project set} also holds when $\mathcal{A}$ is any van Hove sequence.

In \Cref{sec:proof_diffraction_spectrum_of_random_cut-and-project_sets}, we prove the following result for $\mu$-almost every cut-and-project set.
\begin{thm}\label{thm. diffraction spectrum of random cut-and-project sets}
Let $d\geq 2$ and $\mathcal{W}\subset \mathbb{R}^m$ be star-shaped with respect to the origin. Then along any van Hove sequence $\cA$ on $\R^d$ the diffraction spectrum $\widehat{\gamma_{\Lambda(\cW,\cL)_{\vis}}}$ for $\mu$-almost every cut-and-project set $\Lambda(\cW,\cL)$ is a translation bounded pure point measure.
\end{thm}

In view of Theorems \ref{thm:DS_of_visible_points_in_CPS} and \ref{thm. diffraction spectrum of random cut-and-project sets}, we raise the following conjecture:

\begin{conjecture}
Let $d\geq 2$ and $\Lambda=\Lambda(\mathcal{W},\mathcal{L})\subset \mathbb{R}^d$ be an irreducible cut-and-project set. Suppose that the density of $\Lambda_{\vis}$ along a van Hove sequence $\mathcal{A}\subset \mathbb{R}^d$ exists and is positive. Then the diffraction spectrum $\widehat{\gamma_{\Lambda_{\vis}}}$ is a translation bounded pure point measure.
\end{conjecture}

\subsection{Ideas of the proofs}
 Our methods rely on establishing a connection between sets of visible points and weak model sets of maximal density. In the case of simultaneous visibility in $\Z^k$, we realize the set $V(S)$ as a weak model set arising from a suitable cut-and-project set whose internal space is a countable product of finite sets over the prime numbers. Rumsey’s density formula naturally appears as the Haar measure of the corresponding window. We then apply a theorem of Baake, Huck, and Strungaru \cite[Theorem~7]{Strungaru1} on the diffraction of weak model sets of maximal density to deduce that the diffraction spectrum is  translation bounded pure point. The coefficients of the diffraction spectrum are computed by explicitly describing the dual lattice and the characters contained in it. We recall the result of Baake et al. in \Cref{Thm_of_Strungaru1}.

For visible points in cut-and-project sets associated with quadratic number fields, the strategy is similar but requires a more delicate construction of the internal space and window. Using results from \cite[Lemma 5.6]{Barak2}, we characterize visibility in terms of arithmetic conditions and exclusion from a contracted copy of the window. This allows us to construct a new cut-and-project scheme whose associated weak model set of maximal density coincides with the set of visible points.
The proof of \Cref{thm. diffraction spectrum of random cut-and-project sets} proceeds by showing that for $\mu$-a.e. 
cut-and-project set $\Lambda$, the set of visible points $\Lambda_{\vis}$ is a weak model set of maximal density. The result then follows from \Cref{Thm_of_Strungaru1}.
In all of our constructions, the corresponding windows are compact and have empty interior.

\begin{remark}
The sets $V_k$ have also been studied as adelic cut-and-project sets (see \cite[\S~7]{el2017spherical}, \cite[\S~9]{Moody}, and \cite{Sing_slides}). Our construction suggests that the sets $\Lambda(\cW,\cL)_{\vis}$ can likewise be viewed as adelic cut-and-project sets. This perspective  provides a natural framework for studying the actions of $\SL_n(\R)$ and $\ASL_n(\R)$ on the sets $\Lambda(\cW,\cL)_{\vis}$. We plan to return to this direction in future work.
\end{remark}

\subsection{Remarks on the results and work in progress}

Let $K$ be a number field of degree $n = [K: \mathbb{Q}]$ whose ring of integers $\mathcal{O}_K$ is a principal ideal domain. Two distinct points $\x=(x_1,\dots,x_k)$ and $\y=(y_1, \dots, y_k)$ in $\mathcal{O}_K^k$ are said to be mutually visible if
\[
\gcd(\x -\y )=\gcd\big((x_1-y_1),\ldots,(x_k-y_k)\big)=\mathcal{O}_K,
\]
where $(x_i-y_i)$ denotes the ideal generated by $x_i-y_i$.
For a finite set $S \subseteq \mathcal{O}_K^k$, denote by $V(S)$ the set of points that are simultaneously visible from every point in $S$. Analogously to the setting of $\mathbb{Z}^k$, this is defined by
\[
V(S)=\{\x  \in \mathcal{O}_K^k : \gcd(\x -\y )=\mathcal{O}_K \text{ for all } \y \in S\}.
\]
In joint work with W. Takeda \cite{Takeda}, the first-named author shows that the density of $V(S)\subset \mathcal{O}_K^k$ exists along a natural van Hove sequence. This result was already known in the case $S=\{\mathbf{0}\}$; see \cite[Theorem~3.7]{Ferraguti}. Using an argument similar to the proof of \cref{part:1:thm:DS_VS} in \Cref{thm:DS_VS}, it is straightforward to verify that $V(S)$ is a weak model set of maximal density. Consequently, by \Cref{Thm_of_Strungaru1}, {$\widehat{\gamma_{V(S)}}$ is a translation bounded pure point measure along the natural van Hove sequence, provided the density of $V(S)$ is positive.

Furthermore, during the preparation of this manuscript, the authors became aware of recent work by H.~Shamir that generalizes \Cref{theorem on error term}. In particular, Shamir avoids the Hammarhjelm condition and instead only requires the real quadratic field $K= \mathbb{Q}(\sqrt{\mathfrak{d}})$ with $\mathfrak{d} \ge 2$ square-free.

\subsection*{Acknowledgments} 

The authors thank Barak Weiss for many helpful discussions. They are grateful to Michael Baake for bringing \cite{Strungaru1} to their attention, for his valuable comments on an earlier version of this manuscript, and for drawing their attention to the work \cite{araujo2026sarnak}. They also thank Faustin Adiceam, Jens Marklof, and Zeev Rudnick for helpful comments and discussions.
The first author was supported by ISF Grant 2860/24. The second author was supported by ISF Grant 2021/24.

\section{Preliminaries}

We recall the main definitions and results needed for the proofs in this section, including  van Hove sequences, diffraction spectra, cut-and-project sets and weak model sets.

\subsection{Van Hove sequences}\label{subsec:van_Hove_sequences}

Let $G$ be a $\sigma$-compact, locally compact Abelian group (LCAG) with Haar measure $\theta_G$. A sequence $\mathcal{A} = \{A_n\}_{n \in I}$ (typically $I=\mathbb{N}$ or $I=(0,\infty)$) of compact subsets of $G$ is an \emph{averaging sequence} if it is nested ($A_n \subseteq A_{n'}$ whenever $n<n'$), $\vol(A_n) < \infty$ for every $n$, and satisfies $\bigcup_{n \in I} A_n = G$.
Such a sequence is called a \emph{van Hove sequence} if $\theta_G(A_n) > 0$ for all sufficiently large $n$ and for every compact set $K \subseteq G$, the $K$-boundary is asymptotically negligible, that is,
\[
\lim_{n \to \infty}
\frac{\theta_G(\partial^K A_n)}{\theta_G(A_n)}
= 0.
\]
The set $\partial^K B$ denotes the $K$-boundary of $B$, defined by
\[
\partial^K B =
\bigl( (K + B) \setminus \interior{B} \bigr)
\cup
\bigl( (-K + \overline{B^{\mathsf{c}}}) \medcap B \bigr),
\]
where $B^{\mathsf{c}}$, $\interior{B}$ and $\overline{B}$ denote the complement, interior and closure of $B \subset G$, respectively (see \cite[p.~125]{Moody1}, \cite[p.~145]{Martin}).

For example, when $G = \R^d$, the sequence of boxes $A_n = [-n,n]^d$ and the sequence of balls centered at the origin with radius $n$, $A_n = B_n(0)$, are van Hove sequences. The existence of van Hove sequences in $\sigma$-compact LCAGs is shown in \cite{Martin}.

The density of a point set $\Lambda \subset G$, calculated along an averaging sequence $\mathcal{A}=(A_n)_{n \in \mathbb{N}}$, is defined by
\begin{equation}\label{def:density_of_a_set}
\theta(\Lambda)
=
\lim_{n \to \infty}
\frac{\card{\Lambda \medcap A_n}}{\theta_G(A_n)},
\end{equation}
provided that the limit exists.

The proof of \Cref{thm. diffraction spectrum of random cut-and-project sets} requires that the density of the cut-and-project sets is calculated along certain sequence of measurable sets called an \emph{unbounded ordered family}. 
Following \cite[p.6]{Barak1}, a sequence of measurable sets $\{\Omega_T\}_{T>0}$ of $G$ is an unbounded ordered family if 
\begin{enumerate}[label=(\roman*).,ref=(\roman*), leftmargin=*]
    \item \label{item:nested} $0 \le T_1 \le T_2 \Rightarrow \Omega_{T_1} \subset \Omega_{T_2}$,
    \item \label{item:finite} For all $T > 0$, $\vol(\Omega_T) <\infty$,
    \item \label{item:increasing} $\vol(\Omega_T) \to \infty$ as $T \to \infty$, and
    \item \label{item:continuity}For all large enough $V>0$ there is $T$ such that $\vol(\Omega_T)=V$.
\end{enumerate}
Let $\cA= \{A_r\}_{r \in I}$ be a van Hove sequence, then it follows that $\cA$ satisfies \ref*{item:nested}, \ref*{item:finite}, and \ref*{item:increasing}. It also satisfies that
\begin{enumerate}[label=(\roman*~b).,ref=(\roman*~b),leftmargin=*,start=4]
    \item \label{item:continuity2} For all large enough $V'>0$ there is $T$ such that $\vol(A_T) \ge V'$.
\end{enumerate}
By \cite[Lemma~1]{schmidt1960}, there exists an unbounded ordered family $\cA' = \{A'_t\}_{t > 0}$ such that $\cA \subset \cA'$.

\subsection{Cut-and-project sets}
\label{subsec:Generalized_cut_and_project_sets}

Let $G$ and $H$ be $\sigma$-compact LCAGs. Denote their Haar measures by $\theta_G$ and $\theta_H$, respectively.
We assume that $H$ is compactly generated. Hence, $H$ is (up to isomorphism) of the form
$\mathbb{R}^d \times \mathbb{Z}^n \times \mathbb{K}$,
for some integers $d,n \geq 0$ and some compact Abelian group $\mathbb{K}$, see \cite[Theorem~9.8]{HewittRoss}.

A cut-and-project scheme $(G,H,\cL)$ consists of the product group $G \times H$, a full-rank lattice $\cL \subset G \times H$ (that is, a discrete, cocompact subgroup), and the natural projections
\[
\pi_{\mathrm{phys}} : G \times H \longrightarrow G,
\qquad
\pi_{\mathrm{int}} : G \times H \longrightarrow H,
\]
such that
\begin{itemize}
    \item[(I)] The restriction $\pi_{\mathrm{phys}}|_{\cL} : \cL \to G$ is injective;
    \item[(D)] The image $\pi_{\mathrm{int}}(\cL)$ is dense in $H$.
\end{itemize}
Condition~(I) allows us to define the $\star$-map
\begin{equation}
\label{eq:star_map}
(\cdot)^\star :
\pi_{\mathrm{phys}}(\cL)
\longrightarrow
H,
\qquad
x
\longmapsto
x^\star
\df
\pi_{\mathrm{int}}
\bigl(
(\pi_{\mathrm{phys}}|_{\cL})^{-1}(x)
\bigr).
\end{equation}

Let $\cW \subset H$ be a non-empty subset (typically assumed to be measurable), called the \emph{window}. A cut-and-project set for the scheme $(G,H,\cL)$ is defined by
\[
\Lambda(\cW,\cL)
\df
\left\{
\pi_{\mathrm{phys}}(y)
:
y \in \cL,
\,
\pi_{\mathrm{int}}(y) \in \cW
\right\}.
\]
Following \cite{Strungaru1}, a cut-and-project set $\Lambda(\cW,\cL)$ is called a \emph{model set} if the window $\cW$ is relatively compact and has non-empty interior, that is,
\[
\interior{\cW} \neq \emptyset.
\]
If $\cW = \overline{\interior{\cW}}$ and if $\cW$ is compact, then the window is called \emph{proper}.
If $\cW$ is relatively compact, has non-empty interior and satisfies
\[
\theta_H(\partial \cW)=0,
\]
then the cut-and-project set $\Lambda(\cW,\cL)$ is a model set, called \emph{regular}. See \cite[\S~2.1]{Barak1} for further discussion of this terminology. The terms \emph{regular} and \emph{irreducible} cut-and-project sets are used interchangeably in the literature; see \cite[\S~2.1]{Barak1} and \cite[Page~6 and Remark~4]{Strungaru1}. Regular cut-and-project sets are Delone sets (see \cite[p.~744]{Huck1}).

Recall that a set $\mathcal{D} \subset G$ is \emph{uniformly discrete} if there exists an open neighborhood $U \subset G$ of the identity such that any of its translates $xU$, where $x \in G$, contains at most one point of $\mathcal{D}$. A set $\mathcal{D} \subset G$ is \emph{relatively dense} if there exists a compact set $K \subset G$ such that $\mathcal{D}K = G$. A \emph{Delone set} $\mathcal{D} \subset G$ is a uniformly discrete and relatively dense set.

The density of an irreducible cut-and-project set $\Lambda(\cW,\cL)$ along a van Hove sequence $\mathcal{A} = (A_n)_{n \in \mathbb{N}}$ exists and is given by
\begin{equation}\label{eq:density_irreducible_CPS}
\theta(\Lambda(\cW,\cL))
=
\lim_{n\to\infty}
\frac{\card{\Lambda(\cW,\cL)\medcap A_n}}{\theta_G(A_n)}
=
\theta(\cL)\cdot\theta_H(\cW),
\end{equation}
see \cite[Theorem~1]{Moody1}, \cite{Schlotmann} and \cite[Lemma~3.2, Remark~3.3]{Huck1}. We refer to \cite[Chapter~II]{Meyer} for further details on cut-and-project sets.

A cut-and-project set $\Lambda(\cW,\cL)$ is called a \emph{weak model set} if the window $\cW$ is relatively compact and satisfies
\[
\theta_H(\overline{\cW}) > 0.
\]
Every weak model set is uniformly discrete (see \cite[Proposition~2.4]{Huck1}). Moreover, if $\cW$ is nowhere dense, then the cut-and-project set $\Lambda(\cW,\cL)$ has arbitrarily large holes (see \cite[Proposition~2.12]{Huck1}) and therefore is not relatively dense.
We have the following inclusions:
\begin{center}
\catboxI{\shortstack{regular/irreducible\\cut-and-project\\sets}}
\;$\subsetneq$\;
\catboxI{\shortstack{model\\sets}}
\;$\subsetneq$\;
\catboxI{\shortstack{weak model\\sets}}
\;$\subsetneq$\;
\catboxI{\shortstack{cut-and-project\\sets}}
\end{center}

\subsection{Diffraction spectrum}
\label{sec:diff_spectrum}
Let $G$ be a $\sigma$-compact LCAG. A measure $\omega$ on $G$ is said to be \emph{translation bounded} if, for every compact set $K \subseteq G$,
\[
\sup_{t \in G} |\omega|(t + K) < \infty,
\]
where $|\omega|$ denotes the total variation measure of $\omega$. Throughout, measures are viewed as linear functionals on the space $C_c(G)$ of continuous functions with compact support.
Given a measure $\mu$ on $G$ and a function $g \in C_c(G)$, we write
\[
\mu(g)
\df
\int_G g \,\rmd\mu.
\]
Define the twist of a complex valued function $g$ by
\[
\widetilde{g}(x)
\df
\overline{g(-x)}
\]
and the twist of a measure $\mu$ by
\begin{equation}\label{eq:twist_measure}
\widetilde{\mu}(g)
\df
\overline{\mu(\widetilde{g})},
\qquad
g \in C_c(G).
\end{equation}

Let $\mu$ and $\nu$ be measures on $G$. If $\mu$ is translation bounded and $\nu$ is finite, then the convolution $\mu \bigast \nu$ is a measure defined by
\[
(\mu \bigast \nu)(f)
\df
\iint_{G \times G}
f(x+y)\,\rmd\mu(x)\,\rmd\nu(y),
\qquad
f \in C_c(G).
\]
See \cite[Lemma~4.9.19]{moody2017almost} and \cite[Theorem~1.2]{argabright1974fourier}.
Recall that the dual group $\widehat{G}$ is the group of homomorphisms (characters)
\[
\chi : G \longrightarrow \mathbb{S}^1,
\]
equipped with pointwise multiplication and the compact-open topology. The group $\widehat{G}$ is again a LCAG.
Let $\mathcal{A}=(A_n)_{n \in \mathbb{N}}$ be a van Hove sequence in $G$. The \emph{autocorrelation} of $\omega$ along the sequence $\mathcal{A}$ is defined, provided the limit exists, as the weak-$\ast$ limit
\[
\gamma_{\omega}
\df
\lim_{n \to \infty}
\frac{
(\omega|_{A_n})
\bigast
(\widetilde{\omega|_{A_n}})
}{
\theta_G(A_n)
}.
\]
The Fourier transform of the autocorrelation, denoted by $\widehat{\gamma_{\omega}}$, is a positive, translation bounded measure on $\widehat{G}$. It is called the \emph{diffraction spectrum} (or \emph{diffraction measure}) of $\omega$ along the van Hove sequence $\mathcal{A}$.
See  \cite[Theorems~4.10.14~and~4.11.7, and Lemma~4.11.3]{moody2017almost} for further details.

\subsubsection{Further results on diffraction spectrum}
We identify a character $\chi(\cdot)=\chi_u(\cdot)$ with its associated pairing $\langle u,\cdot\rangle$, so that in additive notation we have $\chi_u\chi_{u'}=\chi_{u+u'}$.
There is a natural isomorphism
$\widehat{G \times H}
\cong
\widehat{G}\times\widehat{H}$.
To simplify the notation, denote the canonical projections on the dual groups by
\[
\pi_{\mathrm{phys}} :
\widehat{G}\times\widehat{H}
\longrightarrow
\widehat{G},
\qquad
\pi_{\mathrm{int}} :
\widehat{G}\times\widehat{H}
\longrightarrow
\widehat{H}.
\]

Given the lattice $\cL \subset G \times H$, via the $\star$-map defined in \eqref{eq:star_map}, we may write
\[
\cL
=
\left\{
(x,x^\star)
:
x \in \pi_{\mathrm{phys}}(\cL)
\right\}.
\]
Define $\cL^{0} \subset \widehat{G}\times\widehat{H}$ by
\[
\cL^{0}
\df
\left\{
(\chi_u,\chi_v)\in\widehat{G}\times\widehat{H}
:
\chi_u(x)\chi_v(x^\star)=1
\text{ for all } (x,x^\star)\in\cL
\right\}.
\]
The set $\cL^0$ coincides with the annihilator of $\cL$ in $\widehat{G \times H}$, which in turn implies that $\cL^{0}$ is a lattice.
If $(G,H,\cL)$ is a cut-and-project scheme, then the restriction $\pi_{\mathrm{phys}}|_{\cL^0}$ is injective, and the image
$\pi_{\mathrm{int}}(\cL^0)$
is dense in $\widehat{H}$.
This allows us to define the cut-and-project scheme
$(\widehat{G},\widehat{H},\cL^0)$.
Following \eqref{eq:star_map} and by abuse of notation, we define the $\star$-map
\[
(\cdot)^\star :
\pi_{\mathrm{phys}}(\cL^0)
\longrightarrow
\widehat{H},
\qquad
\chi_u
\longmapsto
\chi_u^\star
\df
\chi_v,
\]
where $\chi_u : \widehat{G} \to \mathbb{S}^1$ and $\chi_v : \widehat{H} \to \mathbb{S}^1$ are as in the definition of $\cL^0$, that is,
\[
\chi_u(x)\chi_v(x^\star)=1
\qquad
\text{for all } (x,x^\star)\in\cL.
\]

The \emph{covariogram function} $c_{\cW}$ of a relatively compact Borel set $\cW \subset H$ is the real-valued function defined by
\[
c_{\cW}(x)
=
(I_{\cW}
\bigast
\widetilde{I_{\cW}})(x),
\qquad
x \in H,
\]
where $I_{\cW}$ denotes the indicator function of $\cW$.
The convolution of two integrable functions $f,g \in L^1(H)$ is defined by
\[
(f \bigast g)(x)
\df
\int_H
f(y)\,g(x-y)\,\rmd\theta_H(y).
\]
It is a standard fact that $c_{\cW}$ is bounded and uniformly continuous, since both $I_{\cW}$ and $\widetilde{I_{\cW}}$ are bounded integrable functions on $H$.

\subsection{Weak model sets of maximal density}
Let $\Lambda(\cW,\cL)$ be a cut-and-project set, where $G$ and $H$ are $\sigma$-compact LCAGs and $H$ is compactly generated.
We say that $\Lambda(\cW,\cL)$ is a \emph{weak model set of maximal density} relative to a given van Hove averaging sequence $\mathcal{A}$ if $\Lambda(\cW,\cL)$ is a weak model set (recall that this means $\cW \subset H$ is relatively compact and $\theta_H(\overline{\cW}) > 0$), and if the density of $\Lambda(\cW,\cL)$ relative to $\mathcal{A}$ exists and is given by
\[
\theta\bigl(\Lambda(\cW,\cL)\bigr)
=
\theta(\cL)\cdot\theta_H(\overline{\cW}).
\]

We will use the following result.

\begin{thm}\emph{\cite[Theorem~7]{Strungaru1}}\label{Thm_of_Strungaru1}
Let $\Lambda = \Lambda(\cW,\cL)$ be a cut-and-project set with compact window $\cW \subset H$ and $\theta_H(\cW) > 0$. If $\Lambda$ is a weak model set of maximal density, then the diffraction spectrum $\widehat{\gamma_{\Lambda}}$ is  translation bounded pure point. Moreover, it is given by
\[
\widehat{\gamma_{\Lambda}}
=
\sum_{\chi_\u \in \pi_{\mathrm{phys}}(\cL^0)}
|a(\chi_\u)|^2\,\delta_{\chi_\u},
\]
where
\[
a(\chi_\u)
=
\frac{\theta(\Lambda)}{\theta_H(\cW)}
\,
\widehat{I_{\cW}(-\chi_\u^\star)},
\]
and $\widehat{I_{\cW}}$ is a bounded, continuous function on $\widehat{H}$.
\end{thm}

Translation boundedness of the diffraction spectrum follows from the fact that it is the  Fourier transform of the autocorrelation which is also translation bounded, see \cite[Theorem~2.5]{argabright1974fourier}.

\section{Proofs}\label{sec:proofs}
\subsection{Proof of \texorpdfstring{\cref*{part:1:thm:DS_VS} of \Cref*{thm:DS_VS}}{~}} \label{proof:thm:DS_VS}

Let $S \subseteq \mathbb{Z}^k$ be a finite set, with $k \geq 2$. We show that $V(S)$, the set of points in $\mathbb{Z}^k$ simultaneously visible from every point of $S$, is a weak model set of maximal density. Then \Cref{Thm_of_Strungaru1} implies that the diffraction spectrum $\widehat{\gamma_{\omega_{V(S)}}}$ is a translation bounded pure point measure.

We work with
\[
G = \mathbb{Z}^k, \qquad H = \prod_{p \in \mathcal{P}} \left( \mathbb{Z}/p\mathbb{Z} \right)^k,
\]
where $H$ is equipped with the product topology and the product normalized Haar probability measure $\theta_H$, so that each factor $(\mathbb{Z}/p\mathbb{Z})^k$ carries the normalized counting measure of total mass~1. The groups $G$ and $H$ are $\sigma$-compact LCAGs, and $H$ is compactly generated. Fix the van Hove sequence $\mathcal{A}= \{[-L,L]^k\}_{L\in \mathbb{N}}$.

Define the set $\mathcal{L}$ by
\[
\mathcal{L} \df \{(\x, \iota(\x)) : \x \in \mathbb{Z}^k\} \subset G \times H,
\]
where
\[
\iota(\x) \df (\pi_p(\x))_{p \in \mathcal{P}}
\]
and $\pi_p : \mathbb{Z}^k \to (\mathbb{Z}/p\mathbb{Z})^k$ is the natural projection defined in \eqref{eq:projection_mod_p}.

The set $\mathcal{L}$ is a discrete subgroup of $G \times H$. Moreover, $\iota(\mathbb{Z}^k)$ is dense in $H$ by the Chinese Remainder Theorem. Since the map 
\[
H \longrightarrow (G \times H)/\mathcal{L}, \qquad h \longmapsto (0,h) + \mathcal{L}
\]
is continuous and surjective and $H$ is compact, it follows that $(G \times H)/\mathcal{L}$ is compact. Hence, $\mathcal{L}$ is a cocompact lattice in $G \times H$.

Let $\pi_{\mathrm{phys}} : G \times H \to G$ and $\pi_{\mathrm{int}} : G \times H \to H$ be the projections. Note that $\pi_{\mathrm{phys}}|_{\mathcal{L}}$ is injective and $\pi_{\mathrm{int}}(\mathcal{L})=\iota(\mathbb{Z}^k)$ is dense in $H$. Hence $(G, H, \mathcal{L})$ is a cut-and-project scheme.

For each $p\in \mathcal{P}$, define
\[
\mathcal{W}_p \df (\mathbb{Z}/p\mathbb{Z})^k \setminus \pi_p(S)
\]
and define
\[
\mathcal{W}_S \df \prod_{p \in \mathcal{P}} \mathcal{W}_p \subset H.
\]
Each $\mathcal{W}_p$ is finite and therefore compact. Hence, by Tychonoff’s theorem, $\mathcal{W}_S$ is a compact subset of $H$. By construction, a vector $\x \in \mathbb{Z}^k$ satisfies $\iota(\x) \in \mathcal{W}_S$ if and only if for every prime $p$, the residue class of $\x$ modulo $p$ does not belong to $\pi_p(S)$. This is equivalent to $\x$ being simultaneously visible from all points of $S$. Therefore, the cut-and-project set
\[
\Lambda(\mathcal{W}_S, \mathcal{L}) \df \left\{ \pi_{\mathrm{phys}}(\y) : \y \in \mathcal{L}, \, \pi_{\mathrm{int}}(\y) \in \mathcal{W}_S \right\}
= \left\{ \x \in \mathbb{Z}^k : \iota(\x) \in \mathcal{W}_S \right\}
\]
coincides with $V(S)$.
The normalized measure of $\mathcal{W}_p$ is
\[
\theta_{H,p}(\mathcal{W}_p) = 1 - \frac{s(p)}{p^k},
\]
where $\theta_{H,p}$ is the normalized Haar measure on $(\mathbb{Z}/p\mathbb{Z})^k$. Hence, by the product measure,
\begin{equation}\label{eq. Haar measure of Ws}
 \theta_H(\mathcal{W}_S) = \prod_{p \in \mathcal{P}} \left(1 - \frac{s(p)}{p^k}\right).   
\end{equation}
This product coincides with Rumsey’s density formula \eqref{eq:sp} of the point set $V(S)$. By the assumption that the density of $V(S)$ is positive, we get $\theta_H(\overline{\mathcal{W}_S})=\theta_H(\mathcal{W}_S) > 0$ and $\Lambda(\mathcal{W}_S, \mathcal{L}) = V(S)$.
By \eqref{eq:sp} and \eqref{eq. Haar measure of Ws}, we get
\[
\theta\bigl(\Lambda(\mathcal{W}_S, \mathcal{L})\bigr) = \theta(V(S))= \theta_H(\mathcal{W}_S)= \theta(\mathcal{L})\theta_H(\overline{\mathcal{W}_S}),
\]
since $\theta(\mathcal{L}) = 1$.
Therefore, $V(S)$ is a weak model set of maximal density. Hence, \Cref{Thm_of_Strungaru1} implies that $\widehat{\gamma_{\omega_{V(S)}}}$ is a translation bounded pure point measure.

\subsection{Proof of \texorpdfstring{\cref*{part:2:thm:DS_VS} of \Cref*{thm:DS_VS} }{~}} \label{proof:thm:DS_VS_2}
Recall that  
$$G = \Z^k, \qquad H= \prod_{p\in \mathcal{P}}H_p,\qquad H_p= \left(\mathbb{Z}/ p\mathbb{Z}\right)^{k},$$
$$\mathcal{W}_{S}= \prod_{p\in \mathcal{P}}\mathcal{W}_{p},\qquad \mathcal{W}_p= H_p\setminus \pi_{p}(S),$$
and the lattice $\mathcal{L}$ is given by
\[
\mathcal{L} \df \{(\x, \x^\star) : \x \in \mathbb{Z}^k\} \subset G \times H,
\]
where $
\x^\star \df \iota(\x) = (\pi_p(\x))_{p \in \mathcal{P}}.$  We use the notation $\x_p$ and $\pi_p(\x)$ interchangeably to shorten the expressions.
The dual lattice of $\mathcal{L}$ has been calculated (see for example \cite[p.6]{Luz} and \cite[p.14]{richard2017short}). For self-containment, we will recall the construction. We begin by setting up the notation for the dual groups. Recall that
\[
\widehat{\mathbb{Z}}
=
\left\{m: x \longmapsto
\chi_m(x) \df \exp{2\pi ixm},\forall x\in \Z 
\right\} = \R.
\]
Since $\widehat{G}
=
\widehat{\mathbb{Z}^k}
\cong
(\widehat{\mathbb{Z}})^k$,
it follows that 
\[
\widehat{G}
\cong
\left\{
\m:
\x \longmapsto
\chi_{\m}(\x)
\df
\exp{2\pi i \langle \x,\m\rangle}~\forall \x \in \Z^k
\right\} = \R^k.
\]
For each prime $p \in \mathcal{P}$, the dual group of $\mathbb{Z}/p\mathbb{Z}$ is
\[
\widehat{\mathbb{Z}/p\mathbb{Z}}
=
\left\{n : y \longmapsto
\chi_n(y)
\df
\exp{2\pi i \frac{yn}{p}},~\forall y \in \Z/p\Z
\right\}= \{0,1,\dots,p-1\}= \mathbb{Z}/p\mathbb{Z}.
\]
Hence
\[
\widehat{(\mathbb{Z}/p\mathbb{Z})^k}
\cong
\left\{\n : \y \longmapsto
\chi_{\n}(\y)
\df
\exp{
2\pi i \frac{\langle \y,\n\rangle}{p}}
\right\} = \{0,1,\dots,p-1\}^k= \left(\mathbb{Z}/p\mathbb{Z}\right)^k.
\]
The dual group $\widehat{H}$ is isomorphic to the direct sum
\[
\widehat{H}
\cong
\bigoplus_{p \in \mathcal{P}}
\widehat{(\mathbb{Z}/p\mathbb{Z})^k}.
\]

An element $(\n_{p})_{p\in \mathcal{P}}\in \widehat{H}$  is a sequence of \emph{finitely many} non-zero vectors. Such sequence represents a  character $\chi_{(\n_p)_{p\in \mathcal{P}}}$ defined by
\[(\y_p)_{p\in \mathcal{P}}\longmapsto \chi_{(\n_p)_{p\in \mathcal{P}}}((\y_p)_{p\in \mathcal{P}})\df \prod_{p\in \mathcal{P}}\exp{2\pi i\frac{\langle\y_p, \n_p\rangle}{p}},\qquad \text{for all }~(\y_p)_{p\in \mathcal{P}} \in H.
\]
Note that the character is well defined, since only finitely many $\n_p$'s are non-zero. Hence the dual lattice of $\mathcal{L}$ consists of points $(\textbf{u}, (\n_{p})_{p\in \mathcal{P}} )\in \widehat{G}\times \widehat{H}$ such that
\begin{equation}\label{eq. annihilator condition}
 \exp{2\pi i \langle \x,\u\rangle}\cdot\prod_{p\in \mathcal{P}}\exp{2\pi i\frac{\langle\x_p, \n_p\rangle}{p}}=1,\qquad \text{for all }~(\x,\x^*)\in \mathcal{L}.    
\end{equation}
Because of linearity and the rules of modular arithmetic, we only need to verify this condition for the canonical basis vectors $\x \in \{e_1,\dots,e_k\}$ of $\mathbb{Z}^k$. By writing the vectors in their component forms as $\u = (u_1, \dots, u_k)$ and $\n_p = (n_{1,p}, \dots, n_{k,p})$, \cref{eq. annihilator condition} breaks down into the following coordinate-by-coordinate requirements  
\begin{equation}\label{eq:indiv_require}
    u_j = \frac{m_j}{N_j} - \sum_{p \in \mathcal{P}} \frac{n_{j,p}}{p}, \qquad \text{ for all } j\in{1,\dots,k},
\end{equation}
where $m_j,~N_j$ are integers and $n_{j,p} \in \{0,\dots,p-1\}$. Consequently, if $(\u,(\n_p)_{p\in \mathcal{P}}) \in \mathcal{L}^0$, then $(\u+\textbf{t}, (\n_p)_{p\in \mathcal{P}}) \in \mathcal{L}^0$ for every translation vector $\textbf{t} \in \mathbb{Z}^k$.  Given that only a finite number of $\n_p$'s are non-zero, the individual conditions of \cref{eq:indiv_require} can be algebraically manipulated to share a least common denominator (\textup{lcm}). Combining all these conditions reveal that they are equivalent to 
\[
\pi_{{\rm phys}}(\mathcal{L}^{0})= \left\{\u =
\frac{\textbf{m}}{N} \in \mathbb{Q}^k : N \text{ is square-free, } \m \in \Z^k, \ \mathrm{gcd}(N,m_1,\dots,m_k)=1 \right\}.
\]
Moreover, the $\star$-map for the dual lattice is given by
\[\u =\frac{\m}{N}\longmapsto \u^\star := \left(-\textup{lcm}(N,p)\cdot \frac{\m}{N} \pmod{p}\right)_{p \in \mathcal{P}}.\]
In this case, the entry at $p$ is $-\m \, \pmod{p}$ for all $p\mid N$ and otherwise $0$. Let $\v_p $ denote the expression $ -\textup{lcm}(N,p)\frac{\m}{N} \pmod{p}$ for each $p \in \cP$.

By \Cref{Thm_of_Strungaru1}, and using the fact that 
$\theta(V(S))= \theta_H(\mathcal{W}_S)$, we have that for all $\u
\in \pi_{\mathrm{phys}}(\mathcal{L}^0)$,
$$a(\u)= \frac{\theta(V(S))}{\theta_H(\mathcal{W}_S)}\cdot \widehat{I_{\mathcal{W}_S}}(-\u^{\star})= \widehat{I_{\mathcal{W}_S}}(-\u^{\star}).$$
By the previous discussion about $\widehat{H}$, $\u^\star = {(\v_p)_{p\in \mathcal{P}}}$ acts via 
\[
\prod_{p \in F} \exp{2\pi i \frac{\langle\v_p,\cdot \rangle}{p}}
\]
over a finite set $F\subset \mathcal{P}$. If $F = \emptyset$, the character is constant and equal to $1$ and if $p \not \in F$, then ${\v_p} = \zerovec\bmod{p}$.

Now, $I_{\mathcal{W}_S}= \prod_{p\in \mathcal{P}}I_{\mathcal{W}_p}$, so the Fourier coefficient factors prime by prime, that is
$$\widehat{I_{\mathcal{W}_S}}(-u^\star)= \prod_{p\in \mathcal{P}}\widehat{I_{\mathcal{W}_p}}(-{\v_p}).$$
For each prime $p$,
\begin{equation}
\begin{split}
 \widehat{I_{\mathcal{W}_p}}(-{\v_p})&=\frac{1}{p^k}\sum_{\x\in H_p\setminus \pi_p(S)}\exp{\frac{-2\pi i \langle {\v_p},\x \rangle}{p}}\\
 &= \frac{1}{p^k}\left(\sum_{\x\in H_p}\exp{\frac{-2\pi i \langle \v_p,\x \rangle}{p}}- \sum_{\x\in \pi_p(S)}\exp{\frac{-2\pi i \langle \v_p,\x \rangle}{p}}\right).
\end{split}    
\end{equation}
Now compute the full character sum. If $\v_p=\zerovec\bmod{p}$, then every exponential is 1, so 
$$\sum_{\x\in H_p}\exp{-2\pi i\frac{\langle \zerovec,\x\rangle}{p}}=\#(H_p) = p^k,$$
and therefore
$$\widehat{I_{\mathcal{W}_p}}(\zerovec)= \frac{p^k-s(p)}{p^k}= 1-\frac{s(p)}{p^k}.$$
If ${\v_p} \neq  \zerovec\bmod{p}$,  choose $j$ such that the $j$-th coordinate of $\v_p$ is non-zero, that is $\v_{p,j}\neq 0 \bmod{p}$. Since $\x\in \{0,\ldots, p-1\}^k$, we get 
$$\sum_{\x\in H_p}\exp{-2\pi i \frac{\langle \v_p,\x\rangle}{p}}= \prod_{r=1}^k\left(\sum_{n=0}^{p-1}\exp{-2\pi i\frac{\v_{p,r} n}{p}}\right).$$
In the $j$-th coordinate, the ratio $\exp{-2\pi i\frac{ \v_{p,j}}{p}}\neq 1$ is a non-trivial $p$-th root of unity, that is
$\sum_{n=0}^{p-1}\exp{-2\pi i\frac{ \v_{p,j} n}{p}}=0$.
Hence, the whole product is $0$. Therefore
$$\widehat{I_{\mathcal{W}_p}}(-{\v_p})= -\frac{1}{p^k}\sum_{\x\in \pi_p(S)}\exp{- 2\pi i\frac{\langle \v_p,\x \rangle}{p}},\qquad {\v_p} \neq \zerovec \bmod{p}.$$
So the local formula is
\begin{equation}
\begin{split}
 \widehat{I_{\mathcal{W}_p}}(-{\v_p})&=\begin{cases}
     1-\frac{s(p)}{p^k},\qquad& {\v_p}= \zerovec\bmod{p},\\
     -\frac{1}{p^k}\displaystyle\sum_{\x\in \pi_p(S)}\exp{- 2\pi i\frac{\langle \v_{p},\x \rangle}{p}},\qquad& {\v_p} \neq \zerovec\bmod{p}.
 \end{cases}    
\end{split}    
\end{equation}
Let $\u= \m/N \in \pi_{\mathrm{phys}}(\cL^{0})$. It follows from the definition of $\mathcal{L}^0$ that only prime divisors of $N$ can contribute non-trivially. Therefore, $$\widehat{I_{\mathcal{W}_S}}(-u^\star)= \prod_{p \nmid N}\left(1-\frac{s(p)}{p^k}\right)\cdot \prod_{p \mid N}\left(-\frac{1}{p^k}\sum_{\x\in \pi_p(S)}\exp{- 2\pi i\frac{\langle \v_p,\x\rangle}{p}}\right).$$
Since 
$$\theta_H(\mathcal{W}_S)= \prod_{p\in \mathcal{P}}\left(1-\frac{s(p)}{p^k}\right)
\quad 
\text{and}
\quad
\theta(V(S))= \theta_H(\mathcal{W}_S),$$
we get
\begin{equation}
\begin{split}
\widehat{I_{\mathcal{W}_S}}(-\u^\star)&= \theta_H(\mathcal{W}_S) \prod_{p\mid N}\frac{\left(-\frac{1}{p^k}\displaystyle\sum_{\x\in \pi_p(S)}\exp{- 2\pi i\frac{\langle {\v_p},\x \rangle}{p}}\right)}{\left(1-\frac{s(p)}{p^k}\right)}\\
&=\theta(V(S)) \prod_{\substack{p\in\mathcal{P}\\ p\mid N}}
\frac{\left(-\frac{1}{p^k}\displaystyle\sum_{\x\in \pi_p(S)}\exp{- 2\pi i\frac{\langle \v_p,\x \rangle}{p}}\right)}{\left(1-\frac{s(p)}{p^k}\right)}.
\end{split}    
\end{equation}

\subsection{Proof of \texorpdfstring{\Cref*{thm:DS_of_visible_points_in_CPS}}{~}} \label{proof:thm:DS_of_visible_points_in_CPS}

Let $\Lambda = \Lambda(\mathcal{W}, \mathcal{L})$ be a cut-and-project set as in the statement of \Cref{thm:DS_of_visible_points_in_CPS}. We show that $\Lambda_{\vis}  = \Lambda(\mathcal{W}, \mathcal{L})_{\vis}$ is a weak model set of maximal density; then \Cref{Thm_of_Strungaru1} implies translation bounded pure point diffraction.

Let $\sigma$ be the non-trivial Galois automorphism of $K$. By abuse of notation,  $\sigma(x_1,\dots,x_d)$ means $(\sigma(x_1)\dots,\sigma(x_d))$. It follows from \cite[Lemma~5.6]{Barak2} that
\begin{equation}\label{eq:inclusion_subset}
\Lambda_{\vis}
= \left\{ \x = (x_1, \ldots, x_d) \in \Lambda(\mathcal{W}, \mathcal{L})_\star :
\gcd(x_1, \ldots, x_d) = \mathcal{O}_K, \; \sigma(\x) \notin \frac{1}{\lambda}\mathcal{W} \right\}.
\end{equation}
This representation allows us to construct a new cut-and-project scheme whose associated weak model set is $\Lambda_{\vis}$.
Define
\[
\mathcal{P}_{\mathcal{O}_K} \df \left\{ \p \subset \mathcal{O}_K : \p \text{ is prime ideal } \right\}.
\]
Let $G = \mathbb{R}^d$ be the physical space and let
\[
H \df \mathbb{R}^d \times \prod_{\p \in \mathcal{P}_{\mathcal{O}_K}} (\mathcal{O}_K/\p)^d
\]
be the internal space with the corresponding projection maps $\pi_{\mathrm{phys}}'$ and  $\pi_{\mathrm{int}}'$.
Both $G$ and $H$ are $\sigma$-compact LCAGs. The group $H$ is compactly generated.

Let $\mink$ be the Minkowski-type embedding
\[
\mink : \mathcal{O}_K^d 
{~\lhook\joinrel\longrightarrow~} G \times H,
\qquad
\x \longmapsto \bigl( \x, \sigma(\x), (\pi_{\p}(\x))_{\p \in \mathcal{P}_{\mathcal{O}_K}} \bigr),
\]
where 
\[
\pi_{\p}: \mathcal{O}_K^d\to \left(\mathcal{O}_K/\p\right)^d,\qquad \p\in \mathcal{P}_{\mathcal{O}_K}
\]
is the quotient map.
Define the set $\mathcal{L}' \df \mink(\mathcal{O}_K^d)$.
The restriction $\pi_{\mathrm{phys}}'|_{\mathcal{L}'}:\cL' \to G$ is injective since the map $\sigma(\cdot)$ determines $\x$ uniquely. 
Moreover, the image $\pi_{\mathrm{int}}'(\mathcal{L}')$ is dense in $H$, because the first factor of $H$ is $\R^d$ and the second factor $\prod_{\p \in \mathcal{P}_{\mathcal{O}_K}} (\mathcal{O}_K/\p)^d$ has the product topology.
The quotient $G \times H / \mathcal{L}'$ is isomorphic to $\R^{2d}/\{(\x,\sigma(\x)):\x \in \mathcal{O}_K^d\}$. Thus, $\cL' \subset G \times H$ is a cocompact lattice.

Define
\[
\mathcal{W}_{\infty}' \df \mathcal{W} \setminus \left(\frac{1}{\lambda}\mathcal{W}\right) = \left \{ w \in \mathcal{W} : w \notin \frac{1}{\lambda}\mathcal{W} \right \}
\]
and for each prime ideal $\p \subset \mathcal{O}_K$, define
\[
\mathcal{W}_{\p}' \df \left(\mathcal{O}_K/\p\right)^d \setminus \{\bf0\}.
\]
Finally, define
\[
\mathcal{W}' \df \mathcal{W}_{\infty}' \times \prod_{\p \in \mathcal{P}_{\mathcal{O}_K}} \mathcal{W}_{\p}'.
\]
By construction, $\mathcal{W}' = \overline{\mathcal{W}'}$. The cut-and-project set associated with this window is
\[
\Lambda' = 
\Lambda'(\mathcal{W}', \mathcal{L}') = \left\{ \pi_{\mathrm{phys}}'(l) : l \in \mathcal{L}', \, \pi_{\mathrm{int}}'(l) \in \mathcal{W}' \right\}.
\]
The set $\Lambda'$ coincides with $\Lambda_{\vis}$. Indeed, the condition $\pi_{\mathrm{int}}'(l) \in \mathcal{W}'$ is equivalent to
\begin{enumerate}
    \item $\x \df \pi_{\mathrm{phys}}'(l) \in \mathcal{O}_K^d$ satisfying $\sigma(\x) \in \mathcal{W}_{\infty}'$, that is, $\x \in \Lambda$ and $\sigma(\x) \notin \tfrac{1}{\lambda}\mathcal{W}$ and
    \item for every prime ideal $\p \subset \mathcal{O}_K$, we have $\x \not\equiv 0 \pmod{\p}$ in $(\mathcal{O}_K/\p)^d$, which is equivalent to $\gcd(\x) = \mathcal{O}_K$.
\end{enumerate}

Next, we compute the measure of $\mathcal{W}'$ in $H$. Each finite factor $\mathcal{W}_{\p}'$, with $\p \subset \mathcal{O}_K$ a prime ideal, has normalized measure
\[
\theta_{H,\p}(\mathcal{W}_{\p}') = 1 - \frac{1}{N(\p)^d},
\]
where $\theta_{H,\p}$ is the Haar probability measure on $(\mathcal{O}_K/\p)^d$. Therefore,
\[
\theta_H(\mathcal{W}') = \vol(\mathcal{W}_{\infty}') \cdot \prod_{\p \in \mathcal{P}_{\mathcal{O}_K}} \left(1 - \frac{1}{N(\p)^d}\right)
= \left(1 - \frac{1}{\lambda^d}\right) \vol(\mathcal{W}) \cdot \frac{1}{\zeta_{\mathcal{O}_K}(d)}.
\]
Since $d \geq 2$, we have $\theta_H(\overline{\mathcal{W}'}) > 0$. The cut-and-project set $\Lambda'(\mathcal{W}', \mathcal{L}')$ is a weak model set.  Given a Jordan measurable set $D$, fix the van Hove sequence $\{TD\}_{T>0}$. It follows from \Cref{theorem on error term} that  $\Lambda'(\mathcal{W}', \mathcal{L}')$ is a weak model set of maximal density. By \Cref{Thm_of_Strungaru1}, the diffraction spectrum of $\Lambda_{\vis}$ is a translation bounded pure point measure.

\subsection{Proof of Theorem \ref{thm. diffraction spectrum of random cut-and-project sets}}\label{sec:proof_diffraction_spectrum_of_random_cut-and-project_sets}
Recall that $\cW \subset \R^m$ is a compact star-shaped set with respect to the origin, with $\cW = \overline{\mathrm{int}(\cW)}$ and boundary of measure zero. Combining \cite[Proposition~3.3~and~Theorem~3.4]{Barak2}, there is a subset $\mathscr{A} \subset \Cl(\R^d)$ of $\mu$-measure $1$ such that:
\begin{enumerate}
    \item For every $\Lambda \in \mathscr{A}$ there exists $\cL \in \scrX_n$ such that $\Lambda = \Lambda(\cW,\cL)$.
    \item Denote by $\mathscr{B}\subset \scrX_n$ the subset of lattices $\cL$ such that $\Lambda(\cW,\cL) \in \mathscr{A}$. This set has $m_{\scrX_n}$-measure 1.
    \item For every $\Lambda \in \mathscr{A}$, we have 
    \[
\Lambda(\cW,\cL_{\mathrm{vis}}) = \Lambda(\cW,\cL)_{\mathrm{vis}}
    \]
    and
    \[
\theta\left({\Lambda(\cW,\cL_{\mathrm{vis}})}\right) = \frac{\vol(\cW)}{\zeta(n)}.
    \]
\end{enumerate}

Let $\b_1,\ldots, \b_n$ be a basis for a lattice $\cL \in \mathscr{B}$. Denote
\[
\Lambda(\mathcal{W},\mathcal{L}_{\vis}) \df \{\pi_{\textup{phys}}(y)\, : \, y\in
\mathcal{L}_{\vis}\,, \pi_{\textup{int}}(y) \in \mathcal{W}\},
\] 
where $\mathcal{L}_{\vis}$ is the set of visible points of the lattice $\mathcal{L}$. 
We define the \emph{content} of a $\y\in \mathcal{L}\setminus\{\zerovec\}$ by
\[
{\rm cont}(\y)\df \max \{l\in \mathbb{N}\, :\, \y\in l\mathcal{L}\}.
\]
If $\y$ is expressed in terms of a basis of $\mathcal{L}$, namely $\y= \sum_{i=1}^n y_i \b_i$, then
\[
{\rm cont}(\y)= \gcd(y_1,\ldots,y_n)
\]
(which is therefore independent of the particular basis chosen). 
For consistency and convenience, we define ${\rm cont}(\textbf{0})=\infty$. 
For $k\in \mathbb{Z}^{+}$ and $\y\in \mathcal{L}$, we have
\[
{\rm cont}(k\y)= k\cdot {\rm cont}(\y),
\]
and
\[
k \mathcal{L}_{\vis}= \{\y \in \mathcal{L}\, :\, {\rm cont}(\y)=k\},
\]
see \cite[p.7]{Moody}.
Denote $G= \mathbb{R}^d$. Note that $\mathcal{L}\cong \mathbb{Z}^n$. For each prime $p$, denote 
\[
H_p=\mathcal{L}/ p \mathcal{L},
\]
and 
\[
H= \mathbb{R}^m\times \prod_{p\in \mathcal{P}}H_p.
\]
We consider $H$ with the product topology. Hence, it is compactly generated, and both $G$ and $H$ are $\sigma$-compact LCAGs. Let $\mathcal{A}$ be a van Hove sequence for $G$.
Define
\[
\mathcal{L}' \df \{(\y, (\y_p)_{p\in \mathcal{P}})\, :\, \y \in \mathcal{L}\},
\]
where 
\[
\y_p \df \sum_{i=1}^n (y_i\bmod{p})\b_i.
\]
Clearly, $\mathcal{L}'$ is a lattice in 
\[
G\times H = \R^d \times \R^m \times \prod_{p\in \mathcal{P}}\mathcal{L}/ p \mathcal{L}.
\]
Define 
$\mathcal{W}_p= \left(\mathcal{L}/ p\mathcal{L}\right)\setminus \{\textbf{0}\},$
and
$$\mathcal{W}'= \mathcal{W}\times \prod_{p\in \mathcal{P}}\mathcal{W}_p.$$
Let $\pi_{\mathrm{phys}}' : G \times H \to G$ and $\pi_{\mathrm{int}}' : G \times H \to H$ be the projections. Note that $\pi_{\mathrm{phys}}'|_{\mathcal{L}'}$ is injective and $\pi_{\mathrm{int}}'(\mathcal{L}')$ is dense in $H$. Hence $(G, H, \mathcal{L}')$ is a cut-and-project scheme.
The cut-and-project set determined by $\mathcal{L}'$ and $\mathcal{W}'$ is
\[
\Lambda' = \Lambda(\cW',\cL')=\{\pi_{{\rm phys}}'(\x)\, :\, \x\in \mathcal{L}'\, , \pi_{{\rm int}}'(\x)\in \mathcal{W}'\}.
\]
The set $\Lambda'$ is equal to $\Lambda(\mathcal{W},\mathcal{L}_{\vis})$.
Indeed, let $\x=(\y,(\y_p)_{p\in\cP}) \in \cL'$ be such that $\pi_{\rm phys}'(\x) \in \Lambda'$. This means that $\y_p\neq \zerovec$ for every prime $p$. Hence
$${\rm cont}(\y)=1.$$ 
Thus $\y\in \mathcal{L}_{\vis}$. By the definition of $\Lambda(\cW,\cL)$, we have $\pi_{{\rm int}}(\y)\in \cW$, and hence $\pi_{{\rm int}}'(\y) \in \mathcal{W}'$.

For each $p\in \mathcal{P}$, the normalized measure
$$\theta_{H,p}(\mathcal{W}_p)= 1-\frac{1}{p^n},$$
where $\theta_{H,p}$ is the Haar probability measure on $\mathcal{L}/ p\mathcal{L}$. Therefore,
$$\theta_H(\mathcal{W}')= \vol(\mathcal{W})\cdot \prod_{p\in \mathcal{P}}\left(1-\frac{1}{p^n}\right)= \frac{\vol(\mathcal{W})}{\zeta(n)}.$$
Since $n\geq 2$, we have $\theta_H(\mathcal{W}')>0$. Hence, the cut-and-project set $\Lambda'$ is a weak model set of maximal density. By \Cref{Thm_of_Strungaru1}, and \eqref{eq. density of the random cut-and-project set} for $\mu$-a.e. cut-and-project set $\Lambda(\mathcal{W},\mathcal{L})$ the diffraction spectrum $\widehat{\gamma_{\Lambda(\cW,\cL)_{\vis}}}$ is a translation bounded pure point measure.

\bibliographystyle{amsalpha}
\bibliography{references}

\end{document}